\documentclass[12pt]{article}
\usepackage{amsmath,amsfonts,amssymb,amsthm,amscd}
\title{On genericity and weight in the free group}
\author{ Anand
Pillay\thanks{This work was supported by a Marie Curie Chair EXC 024052 as well as an EPSRC grant EP/F009712/1}\\University of Leeds, England}
\newtheorem{Theorem}{Theorem}[section]

\newtheorem{Definition}[Theorem]{Definition}
\newtheorem{Remark}[Theorem]{Remark}

\newtheorem{Fact}[Theorem]{Fact}
\newtheorem{Example}[Theorem]{Example}

\newcommand{\Q}{\mathbb Q}  
\newcommand{\Z}{\mathbb Z}

\begin{document}
\maketitle

\begin{abstract}
We prove that the  generic type of the (theory of the) free group $F_{n}$ on $n\geq 2$ generators has infinite weight, strengthening the well-known result that these free groups are not superstable. A preliminary result, possibly of independent interest, is that the realizations in $F_{n}$ of the generic type are precisely the primitives.
\end{abstract} 

\section{Introduction and preliminaries}
From the work by Sela \cite{Sela} and Kharlampovich and Myasnikov \cite{KM}, we know that all nonabelian free groups have the same elementary (or first order) theory, which we call $T_{fg}$. Sela \cite{Sela-stable} also recently proved that $T_{fg}$ is stable. It has been known for a long time (\cite{Gibone}, \cite{Poizat}) that $T_{fg}$ is not superstable.  Rather recently, influenced by work on theories without the independence property, the notion of ``strongly stable" has aroused interest. A theory $T$ is said to be {\em strongly stable} if it is stable and every (finitary) type has finite weight (see below). Any superstable theory is strongly stable. It is natural to ask whether the free group is strongly stable. We show here that it is not. In fact we prove that in $T_{fg}$ the {\em generic type} has infinite weight, strengthening our observation in \cite{Pillay-freegroup} that the generic type has weight $\geq 2$. On the way to proving this we will show that the realizations in $F_{n}$ of the generic type are precisely the primitives. Our proof makes use of recent work by Perin \cite{Perin}.

Many thanks to both Sasha Ivanov and Zlil Sela for helpful discussions and pointing out relevant results in the literature. In particular Ivanov explained to me how the special case ``weight of generic type is $\geq 3$" followed from our methods just using a result of Nies \cite{Nies} that $F_{2}$ is homogeneous.

In the remainder of this section, I recall pertinent facts about free groups, stable groups, and weight. The main results, and their proofs (quite easy) appear in section 2.

\subsection{Free groups}
$F_{n}$ denotes the free group on $n$ generators, where usually $n\geq 2$. In so far as we consider $F_{n}$ as a structure in the sense of model theory it will be in the language of groups $\{\cdot, ^{-1}, e\}$. Sometimes we write $ab$ in place of $a\cdot b$.

It is well-known that $F_{n}$ is isomorphic to $F_{m}$ iff $m=n$. 
If $F_{n}$ is free on generators $a_{1},..,a_{n}$ we call $\{a_{1},..,a_{n}\}$ a basis of $F_{n}$, and an element of a basis is called a {\em primitive} of $F_{n}$. So the primitives of $F_{n}$ form a single orbit under $Aut(F_{n})$.

We will be making heavy use of the following fact, which follows easily from Whitehead's theorem, although we give an explanation.
\begin{Fact} Let $\{a_{1},..,a_{n}\}$ be a basis of $F_{n}$. Let $m\leq n$ and let $k_{1},..,k_{m}$ be integers 
$> 1$. Then $a_{1}^{k_{1}}\cdot .. \cdot a_{m}^{k_{m}}$ is not a primitive of $F_{n}$.
\end{Fact}
\noindent
{\em Explanation.} Let $A = \{a_{1},..,a_{n}\}$. $A^{-1}$ denotes the set of inverses of $A$. An automorphism $\alpha$ of $F_{n}$ is called a Whitehead automorphism if it is induced by either (i) a permutation of $A\cup A^{-1}$, or (ii) for some $i$, a map which fixes $a_{i}$, and for each $j\neq i$ takes $a_{j}$ to $a_{j}a_{i}$, 
$a_{i}^{-1}a_{j}$, or $a_{i}^{-1}a_{j}\cdot a_{i}$. 
A word $w$ (in the basis elements) is said to be cyclically reduced if it is reduced, and not of the form $cvc^{-1}$.
A cyclic word is a cyclically reduced word, defined up to cyclic permutation. For $w$ a reduced word $l(w)$ denotes its length. Whitehead's theorem, which is Theorem 4 in section 10 of \cite{Cohen} or Proposition 4.17 of 
\cite{Lyndon-Schupp} says that if $w,u$ are cyclic words which are in the same orbit of $Aut(F_{n})$ and such that moreover  $l(u)$ is minimal among lengths of words in this orbit, THEN there is a sequence $\tau_{1},..,\tau_{s}$ of Whitehead automorphisms of $F_{n}$ such that for each $i = 1,..,s$, $l(w) \geq l(\tau_{i}...\tau_{1}w)$. Moreover if $l(w) \neq l(u)$, then each of these inequalities is strict.

Now let $w$ be our word $a_{1}^{k_{1}}.....a_{m}^{k_{m}}$. Clearly $w$ is cyclically reduced and so is a  cyclic word. We want to show that $w$ is not primitive, namely is not in the same orbit under $Aut(F_{n})$ as 
$a_{1}$. Supposing otherwise, put $u = a_{1}$, also a cyclic word, clearly of minimal length in its orbit. So Whitehead's theorem applies to $w$ and $u$. But by inspection, $l(\alpha w) \geq l(w)$ for any Whitehead automorphism of $F_{n}$, which gives a contradiction.

\vspace{5mm}
\noindent
In \cite{Sela}, Sela proved:
\begin{Theorem}
If $m \geq n \geq 2$ then the natural embedding of $F_{n}$ in $F_{m}$ is an elementary embedding. 
\end{Theorem}

Of course, this not only solves Tarski's problem on the elementary equivalence of of finite rank nonabelian free groups, but also shows that the natural embeddings of $F_{\kappa}$ into $F_{\lambda}$  (for any finite or infinite cardinals $\kappa < \lambda$) are also elementary embeddings.

On the other hand Chlo\'e Perin \cite{Perin} has recently proved a converse to Theorem 1.2: 
\begin{Theorem} Let $n \geq 2$ and let $G$ be an elementary substructure of $F_{n}$. Then $G$ is a free factor of $F_{n}$. Namely $G$ is a free group of rank at most $n$ and some (any) basis of $G$ extends to a basis of $F_{n}$.
\end{Theorem}

\vspace{2mm}
\noindent
As mentioned in the introduction we let $T_{fg}$ denote common theory of nonabelian free groups of finite rank, which we know now to be complete.
Sela \cite{Sela-stable} has also proved the striking result:

\begin{Theorem} $T_{fg}$ is stable. 
\end{Theorem}

Some basic stability-theoretic properties of the free group will be discussed in the next section.

\subsection{Stability}
For basic model theory we refer the reader to \cite{Marker}. In my recent paper \cite{Pillay-freegroup} on stability-theoretic aspects of the free group, I gave a brief survey of stability and stable groups, directed towards nonexperts, so rather than repeat myself I will direct readers to the introduction of that paper. Further facts about stability and stable groups can be found in \cite{Pillay-book}.

Let us fix a complete countable stable theory $T$, and a model $M$. We let ${\bar M}$ be a saturated elementary extension of $M$ (so we allow $M = {\bar M}$). From stability theory  we have the notion ``$a$ is independent from $b$ over $C$" for $a$, $b$ tuples from $M$ (or even subsets of $M$) and $C$ a subset of $M$. Technically the notion is ``$tp(a/C\cup\{b\})$ does not fork over $C$", where forking is as defined by Shelah. In any case ``$a$ is independent from $b$ over $C$" is synonymous with ``$a$ is independent from $C\cup\{b\}$ over $C$". We also may say 
``$a$ forks with $b$ over $C$" in place of ``$a$ is not independent from $b$ over $C$".

Among key properties of independence in a stable theory are 
\newline
(a) (symmetry) $a$ is independent from $b$ over $C$ iff $b$ is independent from $a$ over $C$,
\newline
(b) (transitivity)  If $A\subseteq B \subseteq C$ then $a$ is independent from $C$ over $A$ iff $a$ is independent from $C$ over $B$ and $a$ is independent from $B$ over $A$. 
\newline
(c) (local character) If $a$ is a finite tuple and $C$ any set then there is a countable $C_{0}\subseteq C$ such that $a$ is independent from $C$ over $C_{0}$.
\newline
(d) (invariance) Whether or not $a$ is independent from $b$ over $C$ depends on $tp(a,b,C)$ (so working in ${\bar M}$ the notion is invariant under automorphism),
\newline 
(e)  (existence) Given $a$ and $C$ and $B\supseteq C$, there is $a'\in {\bar M}$ such that $tp(a'/C) = tp(a/C)$ and $a'$ is independent from $B$ over $C$. 
\newline
(f) (uniqueness)  If $C$ is ``algebraically closed in $M^{eq}$", and $a$ a tuple, then $tp(a/C)$ is {\em stationary}, meaning that for any $B\supseteq C$, if $a', a''$ are such that $tp(a'/C) = tp(a''/C) = tp(a/C)$ and each of $a',a''$ is independent from $B$ over $C$ then $tp(a'/B) = tp(a''/B)$.

\begin{Remark} It would be interesting to give a relatively explicit description of independence when $M$ is a free group. For $C = \emptyset$, and $a,b$ ``generics" this will be done below.
\end{Remark}

\vspace{2mm}
\noindent
Let $A\subseteq M$ and $\{b_{i}:i\in I\}$ be a set of tuples  from $M$. We will say that
$\{b_{i}:i\in I\}$ is $A$-independent, if for each $i\in I$, $b_{i}$ is independent from $\{b_{j}:j\neq i\}$ over $A$.
Let us note a couple of things. Firstly, if $\kappa$ is a cardinal and $(b_{i}:i<\kappa)$ has the property that for each $i$, $b_{i}$ is independent from $\{b_{j}:j<i\}$ over $A$, then in fact $\{b_{i}:i<\kappa\}$ is $A$-independent
(by symmetry and transitivity of independence).

Secondly, if $p(x)\in S(A)$ is stationary, then any two $A$-independent sequences $(b_{i}:i<\kappa)$, $(c_{i}:i<\kappa)$ of realizations of $p$, have the same type over $A$. In particular $(b_{i}:i<\kappa)$ is totally indiscernible over $A$ (every permutation is an elementary map over $A$ in the sense of the ambient model $M$) 

\vspace{5mm}
\noindent
I will give now give a definition of weight in a countable stable theory $T$. Strictly speaking it is the definition of preweight, but for our purposes this will not matter.
\begin{Definition} Let $a$ be a finite tuple from $M$ and $A$ a subset of $M$. The {\em weight} of $tp(a/A)$
(written $w(tp(a/A))$ or even $w(a/A)$) is the supremum of the cardinals $\kappa$ such that in ${\bar M}$ there exists
an $A$-independent set $\{b_{i}:i<\kappa\}$ such that $a$ forks with each $b_{i}$ over $A$.
\end{Definition}

\begin{Fact} For any $a,A$ as in Definition 1.6, $w(a/A) \leq \omega$. Moreover, if $w(a/A) = \omega$ then for some $B\supseteq A$, possibly from ${\bar M}$, such that $a$ is independent from $B$ over $A$, the supremum is achieved for $tp(a/B)$, namely there is a $B$-independent set $\{c_{i}:i<\omega\}$ such that $a$ forks with each $c_{i}$ over $B$.
\end{Fact}

\vspace{2mm}
\noindent
In \cite{Pillay-book} we called a stable theory {\em thin} if every finitary type (namely type of some finite tuple over some set) has finite weight. Following work of
Shelah \cite{Shelah}, this is now called {\em strongly stable}. Any superstable theory is strongly stable. Among the reasons for the current interest in weight is that Shelah was able to find a solution to the equation 
$x:NIP$ = $strongly stable : stable$, rather than a solution to $x:NIP$ = $superstable : stable$.

\begin{Example} (i) In a strongly minimal theory such as the theory of algebraically closed fields of a fixed characteristic, $w(b/A)$ is the same as the Morley rank of $tp(b/A)$  ($b$ a finite tuple, $A$ a set of parameters).
\newline
(ii) If $p(x)$ is a complete stationary type over $\emptyset$ of weight $1$ in a stable theory $T$, then in any model $M$ of $T$ any two maximal independent sets of realizations of $p$ have the same cardinality.
\newline
(iii) If $p$ is a complete type of $U$-rank $\omega^{\alpha}$ (in a stable theory) then $w(p) = 1$. 
\newline
(iv) In the structure $((\Z/4\Z)^{\omega},+)$ the generic type (see below) has Morley rank $2$ but weight $1$. 
\newline
(v) There are strongly stable but non superstable groups, such as a vector space over $\Q$ equipped with predicates for members of an infinite strictly descending chain of subspaces.
\newline
(vi) The generic type (see below) of a separably closed field $F$ of infinite Ersov invariant is $1$, although $Th(F)$ is nonsuperstable.
\end{Example}

\vspace{5mm}
\noindent
For our purposes a {\em stable group} is a definable group in a stable theory $T$. Namely there are formulas $\phi(x)$ and $\psi(x,y,z)$ such that in some (any) model $M$ of $T$, the set of solutions of $\psi$ is the graph of a group operation on the set of solutions of $\phi$. We will assume these formulas have no parameters. So the free group is a stable group, with $T = T_{fg}$. 
In a stable group, the theory of independence above has an equivariant variation, leading to the theory of generic types. Let us fix again a model $M$ of $T$ and let $G$ be the interpretation of the relevant formulas in $M$. There are two equivalent definitions of a ``generic type" or ``generic element" of $G$: (i) Let $A\subseteq M$.
Then $g\in G$ is a generic element of $G$ over $A$, or $tp(g/A)$ is a generic of $G$,  if for any $A$-definable subset $X$ of $G$ containing $g$, finitely many left translates of $X$ cover $G$, (ii) Again for $A\subseteq M$, $g$ is generic in $G$ over $A$ if (working possibly in ${\bar M}$) whenever $h\in G({\bar M})$ is independent from $g$ over $A$, then $h\cdot g$ is independent from $h$ over $A$. 

It is a fact that if $tp(g/A)$ is generic then $g$ is independent from $A$ over $\emptyset$ and also that if $tp(g/A)$ is generic and $g$ is independent from $B$ over $A$, then $tp(g/B)$ is generic.

\vspace{2mm}
\noindent
A stable group is said to be {\em connected} if it has no proper definable subgroup of finite index. Again it is a basic fact that $G$ is connected if and only if there is a unique generic type of $G$: namely, working possibly in ${\bar M}$, for any set $A$ of parameters there is a unique $tp(g/A)$ with $g$ generic in $G$ over $A$. 
For $G$ connected, we denote by $p_{0}^{G}(x)$ the unique generic type of $G$ over $\emptyset$. Moreover this type will be 
{\em stationary}. In particular for any cardinal (maybe finite) $\kappa$, an independent set $\{a_{i}:i\in I\}$ of realizations of $p_{0}^{G}$ is an indiscernible set. It is somehat interesting to note that  the Whitehead transformations applied to such an independent set, are elementary maps. More precisely:
\begin{Fact} Suppose $G$ is a connected stable group. Let $A = \{a_{i}:i\in I\}$ be an independent set of realizations of $p_{0}^{G}(x)$ in $G$. Let $\pi$ be one of the following maps
\newline
(i) for some permutation $\sigma$ of $I$, $\pi(a_{i}) = a_{\sigma(i)}$ or $a_{\sigma(i)}^{-1}$,
\newline
(ii) for some $i\in I$, $\pi(a_{i}) = a_{i}$ and for every $j\neq i$, $\pi(a_{j})$  is $a_{j}\cdot a_{i}$, 
$a_{i}^{-1}\cdot a_{j}$, or $a_{i}^{-1}\cdot a_{j}\cdot a_{i}$. 
\newline
Then $\pi$ is an elementary map in the sense of $G$. In particular $\{\pi(a_{i}):i<\kappa\}$ is also an independent set of realizations of $p_{0}^{G}$. 
\end{Fact}

The following was observed in \cite{Pillay-freegroup}, and moreorever the results were shown to follow 
from Theorems 1.2 and 1.4, using elementary arguments due to Poizat \cite{Poizat}.
\begin{Fact} (i) The free group is connected.
\newline
(ii) If $F$ is a free group with basis $\{a_{i}:i\in I\}$, then $\{a_{i}:i\in I\}$ is an independent set of realizations of the unique generic type $p_{0}$. 
\end{Fact}

Bearing in mind that any two bases of a free group have the same cardinality,  part (ii) above together with Example 1.7 (ii), might be considered evidence that the generic type of the free group has weight $1$. In fact, in \cite{Pillay-freegroup} we already pointed out that the generic type has weight at least two, simply because a generic in the free group is a product of two nongenerics. In the next section we will prove that in fact the generic type of the free group has infinite weight.

\section{Main results}

As above $T_{fg}$ is the theory of the noncommutative free group, and $p_{0}(x)$ is the generic type of $T_{fg}$ over $\emptyset$.
In fact we will be working entirely in ``standard models" of $T_{fg}$, namely free groups of finite rank. 

Our first result is a kind of converse to Fact 1.10 (ii). 
\begin{Theorem} (i) Let $F$ be a free group of finite rank ($\geq 2$). Then any realization of $p_{0}(x)$ in $F$ is a primitive.
\newline
(ii) Any maximal independent set of realizations of $p_{0}(x)$ in $F$ is a basis of $F$.
\end{Theorem}
\noindent
{\em Proof.} (i) Suppose $F$ is free of rank $n$, with basis $\{a_{1},..,a_{n}\}$, and consider $F$ as a subgroup of $F_{n+1}$ where the latter has basis $\{ a_{1},..,a_{n},a_{n+1}\}$. By Theorem 1.2 $F$ is an elementary substructure of $F_{n+1}$. Let $b$ realize $p_{0}$ in $F$. So $b$ also realizes $p_{0}$ in $F_{n+1}$.

\vspace{2mm}
\noindent
{\em Claim I.}  $b,a_{n+1}$ are independent realizations of $p_{0}$ in $F_{n+1}$.  
\newline
{\em Proof.} Work in $F_{n+1}$. By Fact 1.10(ii), $a_{n+1}$ realizes $p_{0}$, and is moreover independent from $\{a_{1},..,a_{n}\}$ over $\emptyset$. But $b\in dcl(a_{1},..,a_{n})$ so $a_{n+1}$ is independent from $b$ over $\emptyset$. This suffices.

\vspace{2mm}
\noindent
Let $G$ be the subgroup of $F_{n+1}$ generated by $\{b,a_{n+1}\}$.

\vspace{2mm}
\noindent
{\em Claim II.} $(b,a_{n+1})$ has the same type in $G$ as in $F_{n+1}$.
\newline
{\em Proof.} By Fact 1.10 (i), and stationarity of $p_{0}$, the type of $(b,a_{n+1})$ in $F_{n+1}$ is the same as the type of a basis of $F_{2}$ in $F_{2}$. Hence $G$ is free with basis $(b,a_{n+1})$, and Claim II follows.

\vspace{2mm}
\noindent
As $\{b,a_{n+1}\}$ generates $G$, it follows from Claim II that $G$ is an elementary substructure of $F_{n+1}$. By Theorem 1.3, $G$ is a free factor of $F_{n+1}$, whence $\{b,a_{n+1}\}$ extends to a basis $\{b,a_{n+1},c_{1},...,c_{n-1}\}$ of $F_{n+1}$. Let $\phi: F_{n+1}\to F$ be the surjective homomorphism defined by: $\phi$ is the identity on $F$ and $\phi(a_{n+1}) = 1$. So $\phi(b) = b$, and  $\{b,\phi(c_{1}),..,\phi(c_{n-1})\}$ generates $F$. As $F$ is free of  rank $n$, by
Proposition 2.7 of \cite{Lyndon-Schupp} for example, $\{b,\phi(c_{1}),..,\phi(c_{n})\}$ is a basis of $F$. So $b$ is a primitive of $F$.

\vspace{2mm}
\noindent
(ii) Let $I$ be a maximal independent set of realizations of $p_{0}$ in $F$. By part (i) and Fact 1.10(ii), 
$|I| > 1$. As in the proof of part (i), the group $G$ generated by $I$ is an elementary substructure of $F$, and is moreover free on $I$. By Theorem 1.3 again, $G$ is a free factor of $F$ whereby $I$ extends to a basis of $F$. Again using Fact 1.10 and our maximality assumption on $I$, $G = F$, so $I$ is a basis
 of $F$. The proof is complete.
 
 \begin{Remark} By part (ii) of the Proposition, any two maximal independent sets of realizations of $p_{0}(x)$ in a ``free, finite rank" model of $T_{fg}$ have the same cardinality, which is again a kind of ``weight $1$" behaviour.
 \end{Remark}
 
 \begin{Theorem} The generic type $p_{0}(x)$ of $T_{fg}$ has infinite weight.
 \end{Theorem}
 \noindent
 {\em Proof.} For each $n\geq 2$ we will find a realization $g$ of $p_{0}$ in $F_{n}$ and independent realizations $b_{1},..,b_{n}$ of $p_{0}$ in $F_{n}$ such that $g$ depends on (forks with) $b_{n}$ for each $n$.
 We will be systematically using Fact 1.9 to check that certain elements we construct are generic, or even independent sets of generics. In fact we could equally well observe that our ``transformations" are taking bases to bases, and hence by Fact 1.10(ii) to independent sets of generics. In any case we will just say ``by Fact 1.9".
 
 \vspace{2mm}
 \noindent
 Let us fix a basis $\{a_{1},..,a_{n}\}$ of $F_{n}$. 
 Let $g = a_{1}a_{2}^{3}a_{3}^{3}...a_{i}^{3}...a_{n}^{3}$. Then
 \newline
 {\em Claim I.}  $g$ realizes $p_{0}(x)$.
 \newline
 {\em Proof.} By Fact 1.9.
 
 \vspace{2mm}
 \noindent
 Now let $b_{1} = a_{1}$, 
\newline
$b_{2} = a_{1}a_{2}$, 
\newline
and for $i=3,..,n$,
\newline
$b_{i} = a_{1}a_{2}^{3}...a_{i-1}^{3}a_{i}$.
\newline
{\em Claim 2.} $\{b_{1},..,b_{n}\}$ is an independent set of realizations of $p_{0}$.
\newline
{\em Proof.}  Again by Fact 1.9. In fact one sees directly that $\{b_{1},..,b_{n}\}$ is a basis of $F_{n}$.

\vspace{2mm}
\noindent
{\em Claim 3.}  $g$ forks with $b_{i}$ for each $i=1,..,n$.
\newline
{\em Proof.} Suppose for a contradiction that $g$ is independent from $b_{i}$. 
By Fact 1.9,  $b_{i}^{-1}g$ realizes
$p_{0}$, so in particular by Proposition 2.1(i) is a primitive element in $F_{n}$. Suppose first that $i=1$. 
As $b_{1} = a_{1}$,  $b_{1}^{-1}g = a_{2}^{3}...a_{n}^{3}$, but by Fact 1.1 the latter is not primitive, contradiction.
Now suppose $i > 1$. Then $b_{i}^{-1}g = a_{i}^{2}a_{i+1}^{3}...a_{n}^{3}$, also not primitive by Fact 1.1. Again a contradiction.

\vspace{2mm}
\noindent
Claims 1, 2 and 3 show that $w(p_{0})\geq n$. As $n$ was arbitrary $w(p_{0}) = \omega$, proving the theorem.

\begin{Remark} It is not hard to conclude from the proof of Theorem 2.3 plus compactness, that in some  model $G$ of 
$T_{fg}$, there is a realization $g$ of $p_{0}$ and an independent set $\{b_{i}:i<\omega\}$ of realizations of $p_{0}$ such that $g$ forks with each $b_{i}$ over $\emptyset$. Can this happen in $F_{\omega}$?
\end{Remark}

\end{document}